\documentclass[11pt, reqno]{amsart}

 \usepackage{amsmath}
 \usepackage{amssymb}
\usepackage{bm}

\newcommand{\mysection}[1]{\section{#1}
      \setcounter{equation}{0}}

\newcommand{\nlimsup}{\operatornamewithlimits{\overline{lim}}}
\newcommand{\nliminf}{\operatornamewithlimits{\underline{lim}}}

\newtheorem{theorem}{Theorem}[section]
\newtheorem{lemma}[theorem]{Lemma}

\newtheorem{corollary}[theorem]{Corollary} 

\theoremstyle{definition}
\newtheorem{assumption}{Assumption}[section]

\theoremstyle{remark}
\newtheorem{remark}{Remark}[section]

\newtheorem{example}{Example}[section]

\newcommand{\loc}{\text{\rm loc}}
\newcommand{\diam}{{\rm diam}}
\newcommand{\sign}{\text{\rm\,sign}\,}
\newcommand\inside{\rm int}

\newcommand\sfu{{\sf u}}

 \makeatletter 
 \def\dashint{%
 \operatorname%
 {\,\,\text{\bf--}\kern-.98em\DOTSI\intop\ilimits@\!\!}}
\def\dashnorm{\,\,\text{\bf--}\kern-.5em\|}
\def\ninf{\qopname\relax\@empty{inf\phantom{p}\!\!\!}}

 \makeatother

\newcommand\gb{\mathfrak{b}}

\newcommand\bbeta{\text{\raise-.2ex\hbox{$\bm{\beta}$}}}

\newcommand\bR{\mathbb{R}}

\newcommand\bL{\mathbb{L}}

\newcommand\bS{\mathbb{S}}

\newcommand\cL{\mathcal{L}}

\newcommand\dist{{\rm dist}\,}

\newcommand\infsup{\operatornamewithlimits{inf\,\,\,sup}}

\begin{document}

\title[Elliptic equations with Morrey drift]
{Linear and fully nonlinear elliptic equations with Morrey drift}

\author{N.V. Krylov}
 
\email{nkrylov@umn.edu}
\address{127 Vincent Hall, University of Minnesota,
 Minneapolis, MN, 55455}



\begin{abstract}
We present some results concerning
the solvability of linear elliptic
equations in bounded domains with the main
coefficients almost in VMO, the drift and the free terms in
  Morrey classes containing $L_{d}$, and
bounded zeroth order coefficient. We prove
that the second-order derivatives of solutions
are in a local Morrey class containing $W^{2}_{p,\loc}$. Actually, the exposition is given for fully nonlinear equations and encompasses
the above mentioned results, which are new even
if the  main part of the equation is just
the Laplacian.

\end{abstract}

\maketitle

\mysection{Introduction}

In this article we are dealing with fully
nonlinear uniformly nondegenerate elliptic equations
\begin{equation}
                        \label{8.10.1}
H[u](x):=H(u,Du,D^{2}u,x)=0
\end{equation}
in bounded domains in $\bR^{d}$, that
is a Euclidean space of points $x=(x^{1},...,x^{d})$,  where $H(\sfu,x)$
is a function given for $x\in\bR^{d}$ and $\sfu=(\sfu',\sfu'')$,
$$
\sfu'=(\sfu'_{0},\sfu'_{1},...,\sfu'_{d}\big)\in\bR^{d+1},
\quad \sfu''\in\bS,
$$  
and $\bS$ is the set of symmetric $d\times d$-matrices. We assume that the growth of $H$ with respect to $|Dv|$
is controlled by the product of $|Dv|$ and a function from a Morrey class {\em containing\/} $L_{d}$. The case when this function is
in $L_{d}$ is treated earlier in \cite{Kr_20}.
The dependence of $H$ on $x$ is assumed to be of BMO type.
Among other things we prove that there exists
$d_{0}\in(d/2,d)$ such that for any $p\in(d_{0},d)$
the equation with prescribed continuous boundary data
has a solution whose second-order
derivatives are locally in a Morrey class
contained in $L_{p}$.
Our results are new
even if $H$ is linear, when \eqref{8.10.1}
becomes
\begin{equation}
                        \label{8.10.2}
a^{ij}(x)D_{ij}+b^{i}(x)D_{i}u+c(x)u +f(x)=0.
\end{equation}
Naturally, 
$$
D_{i}=\frac{\partial}{\partial x^{i}},\quad
D_{ij}=D_{i}D_{j},\quad Du=(D_{i}u),
\quad D^{2}u=(D_{ij}u).
$$

The interest in solutions with $D^{2}u$
in Morrey classes
has a long history (see, for instance,
\cite{FHS_17} for an extensive list of references and the history of the subject). 
The main motivation under this interest is
that on the account of assuming better
summability properties of the free terms
in $L_{p}$, in terms of belonging to a Morrey 
class,
one obtained better summability of $D^{2}u$
in $L_{p}$, in the same terms, and by using
embedding theorems for Morrey classes 
one gets much better regularity of $Du$ than
using just Sobolev embedding theorems.
Roughly speaking, having already estimates
of $D^{2}u$ in $L_{p}$, one tries to extend them to estimates in a Morrey class contained in $L_{p}$.

Our motivation is quite different. If
 $b=(b^{i})$ in \eqref{8.10.2} is not
in $L_{d}$, generally, there are no a priori
estimates, no matter how smooth $u$ is,
even if $a^{ij}=\delta^{ij}$.
\begin{example}
                        \label{example 8.10.1}
Take a constant $c$ and let $b=-cx/|x|^{2}$, $u(x)=1-|x|^{2}$, then
$\Delta u+b^{i}D_{i}u=2c-2d$, which is zero if 
$c=d$. Also $u=0$ on the boundary of $\{|x|<1\}$ and there is no estimates of $u$ through
the free term and boundary data.

\end{example}

In this example $b\not\in L_{d}$, but it is
in Morrey classes containing $L_{d}$ and we will show that, if $c$ is small enough, the equation
$\Delta u+b^{i}D_{i}u=f$ is solvable even when $\Delta$ is replaced
with  $a^{ij}D_{ij}$, if it is a uniformly elliptic operator with $a^{ij}$ almost in VMO.
In short, we can do with the drift term
with summability property
below $L_{d}$.
This possibility was already exploited in the literature, see, for instance, \cite{CGV_98}
and the references therein, however with not so general $a^{ij}$.

The closest to our results the author could find in the literature are those in \cite{FHS_17}
 and \cite{Li_03} which
contain plenty of information, beyond the scope
of this article, in the case of
linear equations with   
$b$ at least in $L_{d}$. For instance, in \cite{FHS_17}
the power of summability $p$ of $D^{2}u$ can be any number
in $(1,\infty)$. In our results when the linear
equation appears as a particular case of fully
nonlinear equation we have a restricted range of $p$, but $b$ is in a Morrey class containing
  $L_{d}$.

The authors of \cite{BLP_16} treat fully nonlinear case by first showing the solvability
in spaces with Muckenhoupt weights and 
from  them,
by using an elegant observation which first
appeared in the proof of
Theorem 3 in \cite{CF_88}, derive the solvability in Morrey
spaces. This approach can also  work in our case, but requires developing first a theory
similar to \cite{DK_19} of solvability
in $W^{2}_{p}$-spaces with weights and $p<d$.
It is assumed in \cite{DK_19} that $p>d$,
but after having the Aleksandrov type
estimate in  \cite{DK_21} and  \cite{Kr_21a} for 
equation  \eqref{8.10.2}
 with $f\in L_{p}$, $p<d$, and $b$ in a Morrey
class containing $L_{d}$ one could mimic
what was done in \cite{DK_19}.

However, we decided to make a shortcut and separate the
Morrey space theory from the weighted space one.
In \cite{BLP_16} the authors obtain global estimates in case the first order
``coefficients'' are bounded. Our estimates are
only local, however, the author does not see
much difficulties to make them global as well.
By the way, specified to \eqref{8.10.2}
the results of \cite{BLP_16} require $a$
to be uniformly sufficiently close to a continuous function. In our case $a\in VMO$
suffices.

We follow a very natural approach
from \cite{Li_03} based on the idea that the Morrey class estimates
should more or less easily follow from the
$L_{p}$ estimates for equations
without lower order terms and then treating lower order terms as perturbations. It is important to stress that
one cannot treat lower order terms
even in Example \ref{example 8.10.1}
 as perturbations  dealing
with the Sobolev space theory {\em but\/}
it is possible to do so in the Morrey space
theory mainly due to the Adams theorem 3.1 of \cite{Ad_75}.  One of main new facts
which we use is the Aleksandrov type
estimate obtained in \cite{Kr_21a}.

The article is organized as follows.
We present our main results and some examples
in Section \ref{section 8.11.1}.
Section \ref{section 8.11.2} contains
some auxiliary results.
In Section \ref{section 8.11.3} we prove
interior estimates and in the final Section
\ref{section 8.11.4} we prove existence theorems.

We finish the introduction with some notation.
Define
$$
B_{r}(x)=\big\{y\in\bR^{d}:|x-y|<r\big\},\quad B_{r}=B_{r}(0).
$$
For measurable $\Gamma\subset \bR^{d}$
set $|\Gamma|$ to be its Lebesgue measure
and when it makes sense set
$$
f_{\Gamma}=\dashint_{\Gamma}f\,dx=\frac{1}{|\Gamma|}
\int_{\Gamma}f\,dx.
$$

For domains $\Omega\subset\bR^{d}$, $p\in[1,\infty)$, and $\mu\in(0,d/p]$, introduce
Morrey's space $E_{p,\mu}(\Omega)$
as the set of $g\in L_{p}(\Omega)$ such that
\begin{equation}
                             \label{8.11.2}
\|g\|_{E_{p,\mu}(\Omega)}:=
\sup_{\rho<\infty,x\in \Omega}\rho^{\mu}
\dashnorm g I_{\Omega}\|_{L_{p}(B_{ \rho}(x) )} <\infty ,
\end{equation}
where
$$
\dashnorm g\|_{L_{p}(\Gamma)}=\Big(
\dashint_{\Gamma}|g|^{p}\,dx\Big)^{1/p}.
$$
Observe that for bounded $\Omega$
one can restrict $\rho$ in \eqref{8.11.2}
to $\rho\leq\diam(\Omega)$ due to $\mu\leq d/p$. 
Let
$$
E^{2}_{p,\mu}(\Omega)=\{u:u,Du,D^{2}u\in E_{p,\mu}(\Omega)\}
$$
and provide $E^{2}_{p,\mu}(\Omega)$ with an obvious norm.
 
We will often, always tacitly, use the following formulas
in which $u(x)=v(x/R)$:
$$
\dashnorm u\|_{L_{p}(B_{R})}
=\dashnorm v\|_{L_{p}(B_{1})},\quad
\|u\|_{E_{p,\mu}(B_{R})}=R^{\mu}
\|v\|_{E_{p,\mu}(B_{1})}, 
$$
$$
\|Du\|_{E_{p,\mu}(B_{R})}=R^{\mu-1}
\|Dv\|_{E_{p,\mu}(B_{1})},\quad
\|D^{2}u\|_{E_{p,\mu}(B_{R})}=R^{\mu-2}
\|v\|_{E_{p,\mu}(B_{1})}.
$$
\mysection{Main results}
                        \label{section 8.11.1}

Fix some constants $\delta\in(0,1]$, $K_{0}
\in[0,\infty)$, and a measurable function
$K(x)\geq0$ on $\bR^{d}$. Let $\bS_{\delta}$
be the subset of $\bS$ consisting
of matrices whose eigenvalues are between 
$\delta$ and $\delta^{-1}$. In our assumptions
below we also use the parameters $\hat\theta,\theta$, $q$, and $\hat b$, the values of which
will be specified in the statements
of our results.
\begin{assumption}
                                    \label{assumption 10.5.1}
There are  measurable  
functions $F(\sfu,x)=F\big( \sfu'',x\big)$  and $G(\sfu,x)$ 
   such that
$$
H =F +G .
$$
 Furthermore, for all $\sfu''\in\bS,\sfu'\in\bR^{d+1}$, and 
$x\in\bR^{d}$, we have 
\begin{equation}
                                                     \label{1,14,2}
\big|G(\sfu,x)\big|\leq \hat\theta|\sfu''|+  K_{0}|\sfu'_{0}|
+\gb(x)|[\sfu']|+K(x),\quad 
F(0,x)\equiv0,
 \end{equation}
where
$$
[\sfu']:=(\sfu'_{1},...,\sfu'_{d}).
$$
 
\end{assumption}

Recall that   Lipschitz
continuous functions are almost everywhere differentiable 
thanks to the Rademacher theorem.

\begin{assumption}
                                \label{assump1} 
  (i) The function $F$ is Lipschitz continuous with respect to $\sfu''$
 and
$D_{\sfu''}F\in \bS_{\delta}$
at 
all points of differentiability
of $F$ with respect to $\sfu''$.

Moreover, there exists  $R_0\in(0,1]$  
such that, if   
$r\in (0, R_0]$, $z\in \bR^{d}$,   
  then
one can find a {\em convex\/} function $\bar{F} (\sfu'' )=
\bar{F}_{z,r } (\sfu'' )$ (independent
of $x $)  for which

(ii) We have $\bar{F}(0)=0$ and $ D_{\sfu'' }\bar{F} \in\bS_{\delta}$
at all points of differentiability
of $\bar{F}$;
 
(iii) 
For any $\sfu''\in\bS$ with $|\sfu''|=1$,  
 we have    
\begin{equation}
                                                \label{7.30.2}
\dashint_{ B_{r}(z)}\sup_{\tau>0}\tau^{-1}
\big|F\big( \tau \sfu'' ,x\big)-\bar{F}(\tau \sfu'')\big| \,dx\leq \theta.
\end{equation}

\end{assumption}

\begin{assumption}
                      \label{assumption 8.2.1}
For any $\sfu''\in\bS$,
$x\in\bR^{d}$, and $\zeta
\in[0,1]$
\begin{equation}
                             \label{8.1.1}
|F(\zeta \sfu'',x)|\leq \zeta |F(\sfu'',x)|+ K (x).
\end{equation}

\end{assumption}
 
\begin{remark}
                        \label{remark 8.11.2}
 Assumption \ref{assumption 8.2.1} is satisfied, for instance, if $F(\sfu,x)$ is boundedly
inhomogeneous in the sense that
$|(\partial/\partial t)\big(tF((1/t)\sfu'',x)\big)|\leq K(x) $ at all points
of differentiability of $tF((1/t)\sfu'',x)$
with respect to $t>0$. Indeed, in that case,
for $t_0\in[0,1]$  
$$
t_{0}F(\sfu'',x) 
=F(t_{0}\sfu'',x)-\int_{t_{0}}^{1}
(\partial/\partial s)\big(sF((1/s)t_{0}\sfu'',x)\,ds
$$
\end{remark}

To finish the setting, take $d_{0}=d_{0}(d,\delta )
\in(d/2,d)$
from \cite{Kr_20}.
In Assumption
\ref{assumption 4.29.1} $q$
is a   number in $[d_{0},d)$.

\begin{assumption}
                    \label{assumption 4.29.1}
There exists a constant $\hat b>0$ such that
for any ball $B$ of radius $r\leq R_{0}$ we have
\begin{equation}
                                                \label{3.26.20}
\dashnorm \gb\|_{L_{q}(B)}\leq \hat b r^{-1}.
\end{equation}
\end{assumption}
\begin{remark}
                        \label{remark 8.11.1}
A simple argument shows that
Assumption \ref{assumption 4.29.1} is satisfied
with any $q<d$ in Example \ref{example 8.10.1} in which
  $b\not\in L_{d,\loc}$. It is also
satisfied with $\hat b$ as small as one likes
on the account of choosing $R_{0}$ small enough
if $\gb\in L_{d}(\bR^{d})$,
or if $\gb$ is just bounded, since by H\"older's
inequality $\dashnorm \gb\|_{L_{q}(B)}
\leq \dashnorm \gb\|_{L_{d}(B)}= N(d)
\|\gb\|_{L_{d}(B)}r^{-1}$.
\end{remark}

 Here is our first main result about the
interior estimates when $\gb$ is bounded
and $p$ is restricted only from below. In it
and below
$\nu=\nu(\alpha,\mu,d,\delta,p)>1$ is taken from Lemma \ref{lemma 7.16.1}.
\begin{theorem}
                      \label{theorem 8.4.1}
Let $p\in(d_{0},\infty)$, $d/p\geq\mu>\alpha>0$ and 
let Assumption \ref{assumption 10.5.1}
be satisfied with   $\hat\theta(\alpha,\mu,d,\delta,p)$
from Lemma \ref{lemma 8.2.1} and
Assumption \ref{assump1}  
be satisfied with $\theta=
\theta(\alpha,d,\delta,p)$ from Lemma
\ref{lemma 5.2.1}.
Suppose that
  Assumption \ref{assumption 8.2.1}   
 is satisfied and
  that    $\gb\leq K_{0}$. Let $\rho\leq
1/\nu$, $u
\in W^{2}_{p,\loc}(B_{\rho})\cap C(\bar B_{\rho})$, $H[u]=0$ in $B_{\rho}$, and
$K\in E_{p,\mu}(B_{\rho })$. Then
$u\in E^{2}_{p,\mu}(B_{\rho_{0}})$,
where $\rho_{0}=\rho/3$,
and   we have
\begin{equation}
                              \label{8.4.2}
\|D^{2}u\|_{E_{p,\mu }(B_{\rho_{0}})}\leq
N\rho^{-2} 
\| u\|_{E_{p,\mu}(B_{\rho})}+N  
\rho^{\mu-2}\sup_{B_{\rho}}|u|+N
\| K\|_{E_{p,\mu}(B_{\rho})}=:I,
\end{equation}
where the constants $N$ depend only on
$K_{0}$, $R_{0}$,
$\alpha,\mu,d,\delta,p$.

\end{theorem}

If $\gb$ is from a Morrey space rather than
being bounded we need to restrict $p$
from above and $\mu$ from below.

\begin{theorem}
                       \label{theorem 8.5.1}
Let $p\in(d_{0},d)$, $d/p>\mu>1$ and 
suppose that Assumption \ref{assump1}  
is satisfied with $\theta=
\theta(\alpha,d,\delta,p)$ for some $\alpha<\mu$
and
  Assumption \ref{assumption 8.2.1}
is also satisfied. Let $\rho\leq
R_{0}\wedge \nu^{-1}$, $u
\in E^{2}_{p,\loc}(B_{\rho})\cap C(\bar B_{\rho})$ and $H[u]=0$ in $B_{\rho}$. Then
there exist $\hat\theta=\hat\theta(\alpha,\mu,d,\delta,p)>0$ and $\hat b=\hat b(\alpha,\mu,d,\delta,p)>0$ such that, if Assumptions
\ref{assumption 10.5.1} and \ref{assumption 4.29.1}
are satisfied with these $\hat\theta$ and 
$\hat b$ and $q=p\mu$ ($\in(d_{0},d)$), 
respectively, then
with  $\rho_{0}=\rho/3$,
    we have
\begin{equation}
                              \label{8.4.20}
\|D^{2}u\|_{E _{p,\mu }(B_{\rho_{0}})}\leq
N\rho^{-2} 
\| u\|_{E_{p,\mu}(B_{\rho})}+N  
\rho^{\mu-2}\sup_{B_{\rho}}|u|+N
\| K\|_{E_{p,\mu}(B_{\rho})} ,
\end{equation}
where the constants $N$ depend only on
  $R_{0}$,
$\alpha,\mu,d,\delta,p$.
\end{theorem}

Observe that Remark \ref{remark 8.11.1} shows that,
for $p\in(d_{0},d)$, $d/p>\mu>1$,
Theorem \ref{theorem 8.5.1} contains
Theorem \ref{theorem 8.4.1}.

The above results are a priori estimates.
Note that in them the equation
$H[u]=0$ is not even assumed to be elliptic.
To state existence theorems we need three more assumptions. 

\begin{assumption}
                     \label{assumption 3.11.2}
We are given a bounded domain $\Omega\subset
\bR^{d}$ such  that for some constants $\rho,\gamma>0$ and 
any $x\in\partial \Omega$ and  $r\in (0,\rho)$ we have $|B_{r}(x)\cap \Omega^{c}|\geq \gamma
|B_{r}|$. We are also given a $g\in C(\partial
\Omega)$.
\end{assumption}

\begin{assumption}
                     \label{assumption 3.11.3}    
 The function   $H(\sfu, x)$ is continuous in $\sfu$
for any $x$, is
  Lipschitz continuous with respect to $\sfu''$, and
$D_{\sfu''}H\in \bS_{\delta}$
at 
all points of differentiability
of $H$ with respect to $\sfu''$.
\end{assumption}

\begin{assumption}
                      \label{assumption 3,7,3}
For all values of the  arguments,
\begin{equation}
                                                     \label{3,7,4}
 H(\sfu',0,x)\sign \sfu'_{0}\leq  
\gb(x)\big|[\sfu']\big|+K(x)\quad (\sign 0:=\pm 1) .
 \end{equation}
 
\end{assumption}

Observe that $H(\sfu,x)$ may not be Lipschitz
continuous with respect to $\sfu'$ or decreasing in $\sfu'_{0}$ in contrast
with conditions in very many articles on the subject,
\cite{BLP_16} including.
We are dealing with the solvability  
of
\begin{equation}
                        \label{7.29.1}
H[u]=0
\end{equation}
in $\Omega$ with boundary condition $u=g$
on $\partial\Omega$.

Here is   our result concerning the solvability
of \eqref{7.29.1} in Morrey  spaces in case $\gb$ is bounded.

\begin{theorem}
                         \label{theorem 10.5.1}
Let $p\in(d_{0},\infty)$, $d/p\geq \mu 
 >0$, and
suppose that $\gb\leq K_{0}$ and Assumptions
\ref{assumption 8.2.1}, \ref{assumption 3.11.2}, \ref{assumption 3.11.3}, and \ref{assumption 3,7,3}
are satisfied. 
 There exist    constants $\hat\theta,
\theta\in (0,1]$, depending    
only on   $d$, $p$, $\delta$,   and  $\mu$,  which
are, generally, smaller than
$\hat\theta,\theta$ from Theorem \ref{theorem 8.4.1} and
such that, if Assumptions \ref{assump1} and \ref{assumption 10.5.1}
are satisfied with these
 $\theta$ and $\hat\theta$, respectively,   
and $K\in E_{p,\mu,\loc}(\Omega)\cap L_{p}(\Omega)$,
then
  there exists  
$u\in E^{ 2}_{p,\loc}(\Omega)\cap C(\bar\Omega)$
 satisfying \eqref{7.29.1} in $\Omega$     and such  
that $u=g $ on $\partial\Omega$. Furthermore, in $\Omega$
\begin{equation}
                                 \label{1.8.2}
|u|   
\leq N\|K\|_{L_{p}(\Omega)}+
\sup_{  \partial \Omega}|g|,
\end{equation}
where $N$ depends only on $p,d,\delta, K_{0}$,
 $R_{0}$, and the diameter of $\Omega$.
\end{theorem}

If $\gb$ is from a Morrey space, the assumptions are stronger.
\begin{theorem}
                         \label{theorem 8.8.1}
                        
Let $p\in(d_{0},d)$, $d/p>\mu>1$,   and 
suppose that  Assumptions
\ref{assumption 8.2.1}, \ref{assumption 3.11.2}, \ref{assumption 3.11.3}, and \ref{assumption 3,7,3}
are satisfied. 
There exist    constants $\hat\theta,
\theta, \hat b\in (0,1]$, depending    
only on   $d$, $p$, $\delta$,   and  $\mu$,  which
are, generally, smaller than
those from Theorem \ref{theorem 8.5.1} and
such that, if Assumptions \ref{assumption 10.5.1}, 
\ref{assump1}, and \ref{assumption 4.29.1}
are satisfied with these
 $\theta$, $\hat\theta$, and $\hat b$, 
$q=p\mu$ ($\in(d_{0},d)$), respectively,   
and $K\in E_{p,\mu,\loc}(\Omega)\cap L_{p}(\Omega)$,
then
  there exists  
$u\in E^{ 2}_{p,\loc}(\Omega)\cap C(\bar\Omega)$
 satisfying \eqref{7.29.1} in $\Omega$     and such  
that $u=g $ on $\partial\Omega$. Furthermore, in $\Omega$
\begin{equation}
                                \label{8.8.3}
|u|   
\leq N\|K\|_{L_{p}(\Omega)}+
\sup_{  \partial \Omega}|g|,
\end{equation}
where $N$ depends only on $p,d,\delta,   R_{0} $,
and the diameter of $\Omega$.
\end{theorem}
\begin{remark}
                        \label{remark 8.12.1}
The fact that we need Assumption \ref{assumption 4.29.1} satisfied with $\hat b$
small enough for \eqref{8.8.3} to hold is well
illustrated by Example \ref{example 8.10.1}.
\end{remark}

\begin{remark}
                        \label{remark 10.30.1}   
Observe that
generally there is no  uniqueness  
in Theorems \ref{theorem 10.5.1} or
\ref{theorem 8.8.1}.
For instance, in the one-dimensional case
the (quasilinear) equation 
$$
 D^{2}u +\sqrt{12|Du |} =0
$$
 for $x\in(-1,1)$
with zero boundary data has two solutions:
one is identically equal to zero and the other one is
 $1-|x|^{3}$.

Another example is given by the (semilinear) equation 
$$
D^{2}u+2u(1+\sin^{2}x+u^{2})^{-1}=0
$$
on $(-\pi/2,\pi/2)$ with zero boundary condition. Again there
are two solutions:
one is $\cos x$ and the other one is identically equal to zero.

\end{remark}

\begin{example}
                                 \label{example 5.16.1}
Let $d=3$, $f,K,\gb\in L_{d}(\Omega)$,
$\alpha\in(0,1]$. Let $w(t)$, $t\in[0,\infty)$,
be a  continuously differentiable
 function  
with sufficiently small derivative.
Then  the   
equation
$$
H(Du,D^{2}u,x):= K(x)\wedge|D_{12}u|
+K(x)\wedge|D_{23}u|+K(x)\wedge|D_{31}u|
 $$
\begin{equation}
                                                      \label{3.4.1}
+ 2\Delta u+w(|D^{2}u|)+\gb(x)|Du|^{\alpha}-f(x)=0
\end{equation}
satisfies our assumptions and
Theorem \ref{theorem 8.8.1} is applicable
with any $p\in(d_{0},d)$, $d/p>\mu>1$.

 Observe that $H$ in \eqref{3.4.1}
is neither convex nor concave with respect to $D^{2}u$.
Also note that we can replace $\Delta u$ with $a^{ij}(x)D_{ij}u$
if $a(x)=(a^{ij}(x))$ is an $\bS_{\delta}$-valued
VMO-function such that $a(x)\geq 2(\delta^{ij})$.
\end{example}
\begin{example}
                                             \label{example 1,13,1}
Let $A$ and $B$ be some countable sets and
assume that for $\alpha\in A$, $\beta\in B$,   $x\in\bR^{d}$,
  and $\sfu'\in\bR^{d+1}$
we are given an $\bS_{\delta}$-valued function 
 $a^{\alpha }( x)$ (independent of $\beta$)
and a real-valued function $b^{\alpha\beta}(\sfu',x)$.
 Assume that 
  these functions are measurable
in $x$, $a^{\alpha }$ and $b^{\alpha\beta}$ are
 continuous with respect to $\sfu' $  
uniformly with respect to $\alpha,\beta,x$,
and  
$$
\big|b^{\alpha\beta}(\sfu',x)\big|\leq \gb(x)
\big|\big(\sfu'_{1},...,\sfu'_{d}\big)\big|+K(x),
$$
where $K,\gb\in L_{d}(\Omega)$.
Next assume that there is an $R_{0}\in(0,\infty)$ such that
for any   $z\in\Omega$, $r\in(0,R_{0}]$  one 
can find
$\bar{a}^{\alpha}\in\bS_{\delta}$ (independent of $x$) such that  
$$
 \dashint_{B_{r}(z)}\sup_{\alpha\in A} 
\big|a^{\alpha} ( x )-\bar{a}^{\alpha}\big|\,dx\leq\theta,
$$
where $\theta$ is sufficiently small
(to accommodate Theorem \ref{theorem 10.5.1}).

Consider equation \eqref{7.29.1},
where
$$
H(\sfu,x):=
\infsup_{\beta\in B\,\,\alpha\in A}
\Big[\sum_{i,j=1}^{d}
a^{\alpha }_{ij}\big( x\big)\sfu''_{ij}+
b^{\alpha\beta}(\sfu',x)\Big].
$$
Define
$$
F(\sfu'',x):=
\sup_{ \alpha\in A}
 \sum_{i,j=1}^{d}
a^{\alpha }_{ij}\big( x\big)\sfu''_{ij}.
$$
As in Example 10.1.24 of \cite{Kr_18} one easily sees that we are in the framework of Theorem \ref{theorem 8.8.1}  
with any $p\in(d_{0},d)$, $d/p>\mu>1$.

\end{example}

\begin{example}
                                              \label{example 1.21.1}
A further specification of Example \ref{example 1,13,1}
is given by linear equations. Suppose that we are given an $\bS_{\delta}$-valued  measurable function
 $a ( x)$  
and an $\bR^{d}$-valued function $b(x)$ such that $b\in E_{p\mu,\mu'}(\bR^{d})$, where
$p\in(d_{0},d)$, $d/p>\mu>1$, $0<\mu'<1$.
 
Next assume that there is an $R_{0}\in(0,\infty)$ such that
for any   ball $B\subset \bR^{d}$ of radius smaller than $R_{0}$
$$
 \dashint_{B }  
 |a(x  )-\bar a_{B} |\,dx\leq\theta.
$$
For $f\in L_{d}(\Omega)$ and
  nonnegative and bounded $c$
 consider   equation \eqref{8.10.2}
in $\Omega$ with boundary condition $u=g$ on 
$\partial\Omega$.

In this situation one can obviously take $F(\sfu'',x)=a^{ij}(x)\sfu''_{ij}$
and satisfy Assumption \ref{assump1}  with $\bar F(\sfu)=
\bar a^{ij}_{B_{r}(z)} \sfu''_{ij}$.
Assumptions \ref{assumption 10.5.1} (with $\hat \theta=0$,
$K_{0}=\sup c$,
$\gb=|b|$, $K=|f|$)  and \ref{assumption 3,7,3}
are also satisfied. 
Then observe that  
$$
r\dashnorm b\|_{L_{p\mu}(B_{r}(x))}
\leq r^{1-\mu'}\|b\|_{E_{p\mu,\mu'}(\bR^{d})},
$$
which for small $r$ can be made as small as we like because $\mu'<1$. Hence,
Assumption \ref{assumption 4.29.1} is satisfied
for small $R_{0}$.
Therefore, by Theorem \ref{theorem 8.8.1},
if $\theta$ is sufficiently small, depending only on
$d,p, \delta$, and $\Omega,g$ satisfy Assumption \ref{assumption 3.11.2}, 
the above boundary value problem has a solution
in $u\in E^{2}_{p,\mu,\loc}(\Omega)\cap C(\bar \Omega)$.
Owing to Theorem \ref{theorem 4.1.1} this solution is
unique. 

Even this result is new.  Also observe that, generally,
$u\not\in W^{2}_{d,\loc}(\Omega)$.
The main novelty in this example is that,
generally, $b
\not\in L_{d,\loc}(\Omega)$,
and even if $a$ is constant the result was not known before.

\end{example}

\mysection{Auxiliary results}
                        \label{section 8.11.2}
Let $\Omega$ be a bounded domain in $\bR^{d}$.
For $\rho>0$ set $\Omega^{\rho}=\{x\in\Omega:
\dist (x,\partial\Omega)>\rho\}$.
\begin{theorem}
                                         \label{theorem 7.22.2}
Let $p>d_{0}$  and $f\in L_{p}(\Omega)$. Then
there exist a   constant  $ \theta\in (0,1]$, depending
only on $d$, $p$, and $\delta$,  
such that, if Assumptions \ref{assump1} is satisfied with this 
 $\theta$, then, for any
 $u\in W^{  2}_{ p ,\loc}(\Omega )\cap C(\bar\Omega )$
 that satisfies $F[u]=f$ in $\Omega $  
and    $0< \rho<\rho_{\inside}(\Omega) \wedge 1$, 
where $\rho_{\inside}(\Omega)$
is the interior radius of $\Omega$,
 we have    
\begin{equation}
                                               \label{7.21.1}
\|u\|_{W^{ 2}_p(\Omega^{\rho})}\le N
\| f\|_{L_p(\Omega )}+N\rho^{-2}\|u\|_{C( \Omega ) },
\end{equation}
where the constants $N$  depend only on
   $d$, $p$, $\delta$,    
$R_{0}$, and
  $ \diam(\Omega)$.      
\end{theorem}

This theorem looks like a particular case of Theorem 1.1 of \cite{Kr_20}, in  which
lower order terms are present in the equation,
but $p$ is restricted to $(d_{0},d)$.
The upper bound $d$ for $p$ is caused by the
presence of the first order terms with
the ``coefficient'' in $L_{d}$. However,
if there are no lower order terms, the arguments
in \cite{Kr_20} go through for the full range
$p>d_{0}$.

Next, introduce $\bL$ as the set of operators
$L=a^{ij}D_{ij}+b^{i}D_{i}$ with measurable
coefficients on $\bR^{d}$ such that $(a^{ij})$
is $\bS_{\delta}$-valued and $|b|\leq\gb$.
Here is Theorem 1.4 of \cite{Kr_21a}.

\begin{theorem}
                          \label{theorem 4.1.1}
Assume that, for a  $R_{0}\in(0,\infty)$,
estimate \eqref{3.26.20} holds with $q=d_{0}$ and
$\hat b=\hat b(d,\delta)$ from Theorem 1.1
of \cite{Kr_21a}
for any $r\in(0,R_{0}]$ and ball
$B $ of radius $r$.
Let  
$u\in W^{2}_{d_{0},\loc}(\Omega)\cap C(\bar \Omega)$, $L\in\bL$. 
Take a function $c\geq0$. Then on $\Omega$
\begin{equation}
                                               \label{4.1.5}
u \leq N 
\|I_{\Omega,u>0}(  Lu-cu)_{-}\|_{L_{d_{0}} }
+\sup_{\partial\Omega}u_{+},
\end{equation}
where $N$ depends only on $d,\delta, R_{0}$, and
the diameter of $ \Omega$.

\end{theorem}

With this version of the Aleksandrov estimate
at hand
one can repeat what is done in \cite{Kr_19_1}
in case $b\in L_{d}$ and arrive at
the following result (see Corollary 4.12 in \cite{Kr_19_1})
about the boundary behavior of solutions
of linear equations.

\begin{lemma}
                                       \label{lemma 7.19.1}
Under the assumptions of Theorem \ref{theorem 4.1.1} let  $0\in\partial \Omega$ 
and suppose that for some constants $\rho,\gamma>0$ and 
any $r\in (0,\rho)$ we have $|B_{r}\cap \Omega^{c}|\geq \gamma
|B_{r}|$. 
Let $w(r)$ be a concave continuous function on $[0,\infty)$
such that $w(0)=0$ and $|u(x)-u(0)|\leq
w(|x|)$ for all $x\in \partial \Omega$.
Then there exists $\beta=\beta(d,\delta, \gamma  )>0$
such that, for any $L\in\bL$ and $x\in \Omega$,
 
\begin{equation}
                                               \label{1.6.4}
|u(x)-u(0)|\leq N|x|^{\beta}\|Lu\|_{L_{d_{0}}(\Omega)}
+\omega\big(N|x|^{\beta/2}),
\end{equation}
where   $N$  depends only on $R_{0}$, $ d,\delta, \gamma,\rho$,  and the diameter
of $\Omega$.
\end{lemma}

The details of the proof of this lemma will
be presented elsewhere.

\mysection{Interior estimates} 
                       \label{section 8.11.3}
  
\begin{remark}
                         \label{remark 7.16.1}
If a function $u$ satisfies $F[u]=f$
in $B_{\rho}$ with $\rho\leq 1$, then
$u_{\rho}(x):=\rho^{-2}u(\rho x)$, satisfies
$F_{\rho}[u_{\rho}]=f_{\rho}$ in $B_{1}$, where
$F_{\rho}(\sfu'',x)=F( \sfu'',\rho x)$,
$f_{\rho}(x)=f(\rho x)$. It is important that
if Assumption \ref{assump1}  
is satisfied   for the original $F$,
then it is also satisfied with the same
 $\theta$ and   $R_{0}$ for $F_{\rho}$.

\end{remark}

\begin{lemma}
                                \label{lemma 5.2.1}
Let $p\in(d_{0},\infty)$, $0<4\rho_{1}\leq
\rho_{2}  \leq 1 $,
$\alpha>0$,   $u\in W^{2}_{p,\loc}(B_{\rho_{2} })
\cap C(\bar B_{\rho_{2} })$ and
$f:=F[u]\in L_{p}(B_{\rho_{2} })$. Then
there exists a constant  $\theta=
\theta(\alpha,d,\delta,p)>0$, such
that, if Assumption \ref{assump1}  
is satisfied with this $\theta$, then
\begin{equation} 
                              \label{7.23.10}
\dashnorm D^{2}u\|_{L_{p}(B_{\rho_{1}})}
\leq N_{1} (\rho_{2}/\rho_{1})^{2}\dashnorm f\|_{L_{p}(B_{\rho_{2}})}
+N  (\rho_{2}/\rho_{1})^{\alpha}\rho_{2}^{-2}\sup_{\partial B_{\rho_{2}}}|u-l|,
\end{equation} 
\begin{equation}
                                  \label{7.15.3}
\dashnorm D^{2}u\|_{L_{p}(B_{\rho_{1}})}
\leq N_{1} (\rho_{2}/\rho_{1})^{2}\dashnorm f\|_{L_{p}(B_{\rho_{2}})}
+N_{2} (\rho_{2}/\rho_{1})^{\alpha}\dashnorm D^{2}u\|_{L_{p}(B_{\rho_{2}})}   , 
\end{equation}
where the constants depend  only on $R_{0}$,
$d,\delta,p$, $\alpha$, and $l$ is any  affine function.
\end{lemma}

Proof. Observe that 
to prove   \eqref{7.15.3} it suffices to concentrate on 
$\alpha<1$. In that case take $q>d/\alpha$ and use
Theorem 10.1.14 of \cite{Kr_18} to find a solution 
$  v\in W^{2}_{q,\loc}(B_{\rho_{2}})\cap C(\bar B_{\rho_{2}})$
of the equation $F[  v]=0$ in $B_{\rho_{2}}$ with boundary
value $  v=u$ on $\partial B_{\rho_{2}}$. 
By taking into account Remark \ref{remark 7.16.1}
and using scaling and
  estimate \eqref{7.21.1}    we get that, 
if $\theta$ is chosen appropriately, then
$$
\dashnorm D^{2}v\|_{L_{q}(B_{ \rho_{2}/2})}
\leq N  (\rho_{2}- \rho_{2}/2)^{-2}\sup_{B_{\rho_{2}}}|v-l|
$$
$$
=N  \rho_{2} ^{-2}\sup_{\partial B_{\rho_{2}}}|v-l|
= N  \rho_{2} ^{-2}\sup_{\partial B_{\rho_{2}}}|u-l| ,
$$ 
where $l$ is any affine function and the equalities follow from the fact that
$0=F[v]-F[0]=a^{ij}D_{ij}v$ for certain 
$\bS_{\delta}$-valued $(a^{ij})$.
It follows  by H\"older's inequality that 
$$
\dashnorm D^{2}v\|_{L_{p}(B_{2\rho_{1}})}
\leq \dashnorm D^{2}v\|_{L_{q}(B_{2\rho_{1}})}
$$
$$
\leq  N  \big(( \rho_{2}/2)/(2 \rho_{1})\big ) ^{ d/q}\dashnorm D^{2}v\|_{L_{q}(B_{ \rho_{2}/2})}
\leq N (\rho_{2}/\rho_{1}) ^{ d/q} \rho_{2}^{-2}\sup_{\partial B_{\rho_{2}}}|u-l|.
$$
Since $d/q<\alpha$, we get
\begin{equation}
                              \label{7.16.4}
\dashnorm D^{2}v\|_{L_{p}(B_{2\rho_{1}})}
\leq N  (\rho_{2}/\rho_{1})^{\alpha}\rho_{2}^{-2}\sup_{\partial B_{\rho_{2}}}|u-l|.
\end{equation}

Then a very particular case of 
\eqref{7.21.1} is that
\begin{equation}
                              \label{7.18.1}
\dashnorm D^{2}u\|_{L_{p}(B_{\rho_{1}})}\leq
N\dashnorm f\|_{L_{p}(B_{2\rho_{1}})}+
N \rho_{1}^{-2}\sup_{B_{2\rho_{1}}}|u-l'|,
\end{equation}
where $l'$ is any affine function.
Here $u-l'=w+(v-l')$, where $w=0$ on $\partial
B_{\rho_{2}}$ and $f=F[w+v]=F[w+v] -F[v]=a^{ij}D_{ij}w $ for some
$\bS_{\delta}$-valued $(a^{ij})$. It follows from \cite{Es_93}, \cite{Ca_95}, \cite{Fo_98} or \cite{Kr_21a}
that $|w|\leq N\rho_{2}^{ 2}\dashnorm f\|_{L_{p}(B_{\rho_{2}})}$. Hence the last term in \eqref{7.18.1} is dominated by   
$$
N(\rho_{2}/\rho_{1})^{ 2}\dashnorm f\|_{L_{p}(B_{\rho_{2}})}+
N \rho_{1}^{-2}\sup_{B_{2\rho_{1}}}|v-l'|,
$$
where the second term, for an appropriate choice of $l'$, is estimated by  a constant times
$\dashnorm D^{2}v\|_{L_{p}(B_{2\rho_{1}})}$
owing to the Poincar\'e inequality.
After that, to get \eqref{7.23.10},
it only remains to refer to \eqref{7.16.4}.
Estimate \eqref{7.15.3} is obtained
from \eqref{7.16.4} by the Poincar\'e inequality.
The lemma is proved.

The following result is well known when
$p\geq d$ and $F=\bar F$.
\begin{corollary}
                    \label{corollary 8.4.1}
If $F$ is independent of $x$, $p\in(d_{0},\infty)$, $\alpha>0$, 
Assumption \ref{assump1}  
is satisfied with $\theta(\alpha,d,\delta,p) $,
$u\in W^{2}_{p,\loc}(\bR^{d})$ satisfies
$F[u]=0$ in $\bR^{d}$ and
$$
\nliminf_{R\to\infty}R^{\alpha-2}\inf_{l\in\cL}
\sup_{B_{R}}|u-l|=0,
$$
where $\cL$ is the set of affine functions,
then $u\in \cL$.

\end{corollary}

\begin{assumption}
                     \label{assumption 7.31.1}
We have $p\in(d_{0},\infty)$, $\mu>\alpha>0$ and 
Assumption \ref{assump1}  
is satisfied with $\theta=
\theta(\alpha,d,\delta,p)$ introduced in Lemma \ref{lemma 5.2.1}.
\end{assumption}

\begin{lemma}
                     \label{lemma 7.16.1}
Suppose that
Assumption \ref{assumption 7.31.1}  
is satisfied. Let
  $\rho \leq 1 $,
$u\in W^{2}_{p,\loc}(B_{\rho  })
\cap C(\bar B_{\rho  })$ and
$f:=F[u]\in L_{p}(B_{\rho  })$.
 Then
there exists $\nu=\nu(\alpha,\mu,d,\delta,p)>1$
such that
for any  $r\leq\rho /\nu$
we have
$$
r^{\mu}\dashnorm D^{2}u\|_{L_{p}(B_{r})}
\leq N \sup_{r\leq s\leq\rho }
s^{\mu}\dashnorm f\|_{L_{p}(B_{s })}
$$
\begin{equation}
                                 \label{7.18.2}
+(1/2)\min\big(\rho ^{\mu}\dashnorm D^{2}u\|_{L_{p}(B_{\rho })},N\rho ^{\mu-2}\sup_{B_{\rho }}
|u-l|\big),
\end{equation}
where $l$ is any affine function and the constants
 $N$ depend only on $\alpha$, $\mu$, 
$R_{0}$,
$d,\delta,p$.
\end{lemma}

Proof. Take the smallest  $\kappa\geq 4$ such that
$N_{2}  \kappa^{\alpha-\mu }\leq 1/2$.
Then, for $r\leq \rho/\kappa$, define 
$r_{n}=\kappa^{n}r$, $m=\max\{n\geq 0:
r_{n+1}\leq\rho \}$,
$$
A_{n}=r_{n}^{\mu}\dashnorm D^{2}u\|_{L_{p}(B_{r_{n}})},\quad B =\sup_{r\leq s\leq\rho }
\rho ^{\mu}\dashnorm f\|_{L_{p}(B_{s})}.
$$
For $0\leq n\leq m$, $\rho_{1}=r_{n}$,
$\rho_{2}=r_{n+1}$ estimate \eqref{7.15.3}
yields
$$
A_{n }\leq N_{1}\kappa^{2-\mu}B +N_{2}\kappa^{\alpha-\mu}
A_{n+1}\leq N_{3}B +(1/2)
A_{n+1}.
$$
By iterating we obtain
$A_{0}\leq 2N_{3}B+2^{- m}A_{m }$ and arrive at
\begin{equation}
                              \label{7.22.1}
r^{\mu}\dashnorm D^{2}u\|_{L_{p}(B_{r})}
\leq N \sup_{r\leq s\leq\rho }
s ^{\mu}\dashnorm f\|_{L_{p}(B_{s })}
+2^{- m}r_{m}^{\mu}\dashnorm
D^{2}u\|_{L_{p}(B_{r_{m}})}.
\end{equation}
Here $r_{m}\leq\rho /\kappa$ and
$r_{m}\geq \rho /\kappa^{2}$. It follows that the last term in \eqref{7.22.1} is less than  
$$
2^{ -m}r_{m}^{\mu}(r_{m}/\rho )^{-d/p}
\dashnorm
D^{2}u\|_{L_{p}(B_{\rho })}\leq
2^{ -m}\kappa^{-\mu+2d/p}
\rho ^{\mu}\dashnorm
D^{2}u\|_{L_{p}(B_{\rho })}.
$$
Now define $m_{0}=m_{0}(\mu,d,\delta,p)$
as the smallest integer $m\geq1$ satisfying the inequality
$2^{ -m}\kappa^{-\mu+2d/p}\leq 1/4$ and set
$\nu=\kappa^{m_{0}}$. Then for $r\leq\rho /\nu$ we have $m\geq m_{0}-1$ and 
$2^{ -m}\kappa^{-\mu+2d/p}\leq 1/2$, so that
the left-hand side of \eqref{7.18.2}
is less than the first term on the right plus one half of the first term inside the min sign. On the other hand,
$$
\dashnorm
D^{2}u\|_{L_{p}(B_{r_{m}})}\leq N\dashnorm
f\|_{L_{p}(B_{2r_{m} })}+
Nr_{m}^{-2}\sup_{B_{2r_{m} }}|u-l|,
$$
and the left-hand side of \eqref{7.18.2}
is less than the first term on the right plus one half of the second term inside the min sign as well. This proves the lemma. 

Below by $\nu$ we always mean the constant from 
Lemma \ref{lemma 7.16.1}.

The following is quite natural. It looks like
it first appeared in \cite{FHS_17}, 
which makes the author wonder how in the past people
claiming that flat boundary and interior
Morrey estimates lead to global estimates
in smooth domains using flattening the boundary 
and partitions of unity. These procedures
unavoidably lead to appearance of the first
order terms, the way to deal with which was not exhibited before \cite{FHS_17}. Unfortunately,
the proof in \cite{FHS_17} contains  an error 
(see Lemma 4.2 there). We give  a different proof.  
\begin{lemma}
                      \label{lemma 7.30.1}
Let $p\in(1,\infty)$, $0<\mu\leq d/p$,
$R\in(0,\infty)$, $u\in W^{2}_{p}(B_{R})$,
and $x_{0}\in B_{R}$. Then there is a constant
$N=N( d,p,\mu)$ such that, for any $\varepsilon
\in(0,1] $, $r\leq 2R, $
$$
 r^{\mu}
\dashnorm I_{B_{R}}Du\|_{L_{p}(B_{r}(x_{0}))}
\leq 
 N\varepsilon R \sup_{r\leq s\leq 2R}s^{\mu}
\dashnorm I_{B_{ R}}D^{2}u\|_{L_{p}(B_{s}(x_{0}))}
$$
\begin{equation}
                              \label{7.30.1}
+N\varepsilon^{-1}R^{-1}
\sup_{r\leq s\leq 2R}s^{\mu}\dashnorm I_{B_{R}}( u-c)\|_{L_{p}(B_{s}(x_{0}))},
\end{equation}
where $c$ is any constant. In particular,
\begin{equation}
                              \label{8.4.3}
\|Du\|_{E_{p,\mu}(B_{R})}\leq
N\varepsilon R \|D^{2}u\|_{E_{p,\mu}(B_{R})}
+N\varepsilon^{-1}R^{-1}
\| u\|_{E_{p,\mu}(B_{R})}.
\end{equation}
\end{lemma}

Proof. Scalings show that we may assume that $R=1$. Obviously we may also assume that $c=0$.
\newpage

 Then denote $v= Du$,
$w= D^{2}u$, $G_{s}=B_{s}(x_{0}) \cap C_{1} $,
$$
U=\sup_{r\leq s\leq 2 }s^{\mu}
\dashnorm I_{B_{1}}u\|_{L_{p}(B_{s}(x_{0}))},\quad
W=\sup_{r\leq s\leq 2 }s^{\mu}
\dashnorm I_{B_{1}}D^{2}u\|_{L_{p}(B_{s}(x_{0}))},
$$
 
By Poincar\'e's inequality, for $ r\leq s\leq  2$,
$$
\dashnorm v-v_{G_{s}}\|_{L_{p}(G_{s})}\leq
N(d)s\dashnorm w\|_{L_{p}(G_{s})}
\leq Ns^{1-\mu}W.
$$
Also by interpolation inequalities,
   there
exists a constant $N=N( d,p)$
such that, for  
 $\varepsilon \in(0,1]$ and $\varepsilon\leq s\leq  2$ ,
$$
\dashnorm v-v_{G_{s}}\|_{L_{p}(G_{s})}
\leq 2\dashnorm v\|_{L_{p}(G_{s})}
\leq N \dashnorm w\|^{1/2}_{L_{p}(G_{s})}
\dashnorm  u\|^{1/2}_{L_{p}(G_{s})}
$$
\begin{equation}
                                 \label{7.31.1}
+Ns^{-1}\dashnorm  u\| _{L_{p}(G_{s})}
\leq N \dashnorm w\|^{1/2}_{L_{p}(G_{s})}
\dashnorm  u\|^{1/2}_{L_{p}(G_{s})}
+N\varepsilon^{-1}\dashnorm  u\| _{L_{p}(G_{s})},
\end{equation}
which for $2\geq s\geq \varepsilon\vee r$    yields
$$
s^{\mu}\dashnorm v-v_{G_{s}}\|_{L_{p}(G_{s})}
\leq N W^{1/2}U^{1/2}   
+N \varepsilon^{-1} U.
$$
  Hence, for any $\varepsilon\in(0,1]$ 
and $r\leq s\leq 2$
$$
s^{\mu}\dashnorm v-v_{G_{s}}\|_{L_{p}(G_{s})}
\leq N_{1}\varepsilon W 
+N_{2}\varepsilon^{-1}U,
$$
where $N_{1}=N_{1}(d,p)$, 
$N_{2}=N_{2}( d,p)$.

Since $\mu\in(0,d/p]$, Campanato's
results (see, for instance, Proposition 5.4 in \cite{MM_12}) imply that  
$$
r^{\mu}\dashnorm v \|_{L_{p}(G_{r})}
\leq N_{3}( N_{1}\varepsilon W 
+N_{2}\varepsilon^{-1}U)+N_{3}\dashnorm v \|_{L_{p}(G_{2})},
$$
where $N_{3}=N_{3}(d,p,\mu)$. 
We estimate the last term as in \eqref{7.31.1}
and come to what  implies 
\eqref{7.30.1}.
The lemma is proved.
\begin{lemma}
                      \label{lemma 8.1.1}
Let $\mu\leq d/p$ and suppose that
  Assumptions \ref{assumption 8.2.1}
and \ref{assumption 7.31.1}
are satisfied and we have
  $r\leq \rho_{0}<\rho_{1}\leq 1/\nu$, 
such that $\rho_{1}-\rho_{0}\leq\rho_{0}$. 
Let $u\in W^{2}_{p}(B_{\rho_{1}})$.
Set $f=F[u]$.
Then 
\begin{equation}
                                 \label{8.1.2}
r^{\mu}\dashnorm D^{2}u\|_{L_{p}(B_{r})}  
\leq N \bar f 
 +N(\rho_{1}-\rho_{0})^{-2} 
\bar u+N\hat u+N\bar K  ,
\end{equation}
where
$$
\bar f:=\sup_{r\leq s\leq\rho_{1} }
s^{\mu}\dashnorm f\|_{L_{p}(B_{s })},\quad
\bar u:=\sup_{r\leq s\leq\rho_{1} }
s^{\mu}\dashnorm u\|_{L_{p}(B_{s })},
$$
$$
\hat u:=\rho_{1}^{\mu-2}\sup_{B_{\rho_{1}}}| u|,\quad
\bar K :=\sup_{r\leq s\leq\rho_{1} }
s^{\mu}\dashnorm K \|_{L_{p}(B_{s })}
$$
and the constants depend only on $\alpha,\mu,
d,p,\delta, R_{0}$.
\end{lemma}

Proof. Define $\kappa=\rho_{1}-\rho_{0}$, $r_{0}=\rho_{0}$ and for $n\geq1$
$$
r_{n}=\rho_{0}+\kappa\sum_{k=1}^{n}
2^{-k}.
$$  
Also introduce smooth $\zeta_{n}(x)$
such that $0\leq \zeta_{n}(x)\leq 1$, 
$\zeta_{n}(x)=1$ for $|x|\leq r_{n}$,
$\zeta_{n}(x)=0$ for $|x|\geq r_{n+1}$,
$$
|D\zeta_{n} |\leq N(d)2^{n}\kappa^{-1},\quad
|D^{2}\zeta_{n} |\leq N(d)4^{n}\kappa^{-2}.
$$

Observe that $r_{n}\le \rho_1= (\nu\rho_{1})/\nu$
and $\nu \rho_{1}\leq 1$, so that 
by Lemma \ref{lemma 7.16.1}
\begin{equation}
                            \label{8.1.5}
W_{n}:=\sup_{r\leq s\leq r_{n}}
s^{\mu}\dashnorm D^{2}u\|_{L_{p}(B_{s})}
\leq N\hat u+NF_{n},
\end{equation}
where
$$
F_{n}:=\sup_{r\leq s\leq \nu\rho_{1}}
s^{\mu}\dashnorm F(D^{2}(\zeta_{n}u))\|_{L_{p}(B_{s})}=
\sup_{r\leq s\leq r_{n+1}}
s^{\mu}\dashnorm F(D^{2}(\zeta_{n}u))\|_{L_{p}(B_{s})},
$$
where the equality is due to H\"older's
inequality, assumption that $\mu\leq d/p$,
and the fact that $F(D^{2}(\zeta_{n}u))=0$
outside $B_{r_{n+1}}$.
By our assumptions (dropping the argument $x$)
$$
|F(D^{2}(\zeta_{n} u) ) |
=|F(\zeta_{n}D^{2}u+2D\zeta_{n}\otimes Du
+uD^{2}\zeta_{n})|
$$
$$
\leq  \zeta_{n}| F(D^{2}u )|+
NI_{B_{r_{n+1}}}(2^{n}\kappa^{-1}|Du|+4^{n}\kappa^{-2}|u|)+K 
I_{B_{r_{n+1}}},
$$
which implies that  
$$
F_{n}\leq N\bar f
+N2^{n}\kappa^{-1}
\sup_{r\leq s\leq r_{n+1}  }
s^{\mu}\dashnorm Du\|_{L_{p}(B_{s })} 
+N4^{n}\kappa^{-2}\bar u
+N\bar K .
$$

Here for $r\leq s\leq r_{n+1}$  by Lemma \ref{lemma 7.30.1} (with $x_{0}=0$, $R=r_{n+1}$)

\begin{equation}
                               \label{8.3.1}
 s^{\mu}
\dashnorm  Du\|_{L_{p}(B_{s} )}
\leq N\varepsilon r_{n+1} W_{n+1}
+N\varepsilon^{-1}r_{n+1}^{-1} \bar u,
\end{equation}
so that for any $\varepsilon\in(0,1]$
$$
F_{n}\leq N\bar f+N_{1}\varepsilon 2^{n}
\kappa^{-1}r_{n+1}
W_{n+1}+N(2^{n}\kappa^{-1}\varepsilon^{-1}r^{-1}_{n+1}
+4^{n}\kappa^{-2} ) 
\bar u+N\bar K .
$$
We may assume that $N_{1}\geq1$, so that
$N_{1}  2^{n}\kappa^{-1}r_{n+1}\geq \rho_{0}
/(\rho_{1}-\rho_{0})\geq1$. Therefore,
now we can take $\varepsilon\in(0,1]$ such that 
$$
N_{1}\varepsilon 2^{n}\kappa^{-1}r_{n+1}= 1/8.
$$
Then $\varepsilon^{-1}= N2^{n}
\kappa^{-1}r_{n+1} $
and
$$
2^{n}\kappa^{-1}\varepsilon^{-1}r^{-1}_{n+1}
+4^{n}\kappa^{-2}\leq N4^{n}\kappa^{-2} .
$$ 
Coming back to \eqref{8.1.5} we get
$$
W_{n}\leq  N(\bar f+\hat u)+8^{-1}
W_{n+1}+N4^{n}\kappa^{-2} 
\bar u+N\bar K .
$$
We multiply this inequality by $8^{-n}$
and sum over $n=0,1,...$. Then we obtain
$$
\sum_{n=0}^{\infty}8^{-n}W_{n}\leq  N(\bar f+\hat u)+\sum_{n=1}^{\infty}8^{-n}
W_{n }+N \kappa^{-2} 
\bar u+N\bar K .
$$
Canceling (finite) like terms yields
$W_{0}\leq  N(\bar f+\hat u)+N \kappa^{-2} 
\bar u+N\bar K $, which implies \eqref{8.1.2}
and proves the lemma.

\begin{lemma}
                      \label{lemma 8.2.1}
Under the assumptions of Lemma \ref{lemma 8.1.1}
suppose that $\gb\leq K_{0}$ and $H[u]=0$ in $B_{\rho_{1}}$. Then there exists $\hat \theta=
\hat\theta(\alpha,\mu,d,\delta,p)>0$ such that,
if Assumption \ref{assumption 10.5.1}
is satisfied with this $\hat\theta$, then
\begin{equation}
                                 \label{8.2.1}
r^{\mu}\dashnorm D^{2}u\|_{L_{p}(B_{r})}
\leq  N (\rho_{1}-\rho_{0})^{-2}  
\bar u+N\hat u+N\bar K  ,
\end{equation}
where
 $N$ depends only on $K_{0}$, $R_{0}$,
$\alpha,\mu,d,\delta,p$.
 
\end{lemma}

Proof. We use the same notation as in   Lemma \ref{lemma 8.1.1} and observe that 
in light of Assumption \ref{assumption 10.5.1},
\begin{equation}
                              \label{8.5.2}
|F(D^{2}u)|
\leq  \hat\theta |D^{2}u|+K_{0}| u|
+  K 
+ \gb |Du|=:f.
\end{equation}
By Lemma \ref{lemma 8.1.1} for any $n\geq0$
and $r\leq r_{n}$
\begin{equation}
                                 \label{8.1.20}
r^{\mu}\dashnorm D^{2}u\|_{L_{p}(B_{r})}
\leq N \bar f _{n+1}
 +N\kappa^{-2} 
\bar u+N\hat u+N\bar K  ,
\end{equation}
where
$$
\bar f_{n+1}:=  \sup_{r\leq s\leq r_{n+1} }
s^{\mu}\dashnorm f\|_{L_{p}(B_{s })}.
$$
Note that
$$
\sup_{r\leq s\leq r_{n+1} }
s^{\mu}\dashnorm \gb |Du|\|_{L_{p}(B_{s })}
\leq K_{0}\sup_{r\leq s\leq r_{n+1} }
s^{\mu}\dashnorm  Du \|_{L_{p}(B_{s })},
$$
where in light of \eqref{8.3.1} the last term
is dominated by
$$ N\varepsilon r_{n+1} W_{n+1}
+N\varepsilon^{-1}r_{n+1}^{-1} \bar u
$$
for any $\varepsilon\in(0,1] $
with $N$ depending only on $ d,p,\mu$.
Hence,
$$
W_{n}\leq N_{1}(\hat\theta+\varepsilon
r_{n+1})W_{n+1}+N(\kappa^{-2}+\varepsilon^{-1}
r_{n+1}^{-1} )\bar u+N\hat u+N \bar K .
$$
We choose here $\hat\theta$ so that $N_{1}
\hat\theta\leq 1/4$ and choose the largest $\varepsilon\in(0,1] $ such that $N_{1}\varepsilon
r_{n+1}\leq 1/4$. Observe that 
$r_{n+1}^{-1}
\leq \kappa^{-1}$, so that in any case
$\varepsilon^{-1}
r_{n+1}^{-1}\leq N+\kappa^{-1}$ and,
since $\rho_{1}<1$, we  have $\kappa\leq1$,
$\kappa^{-1}\leq\kappa^{-2}$, and
 $$
W_{n}\leq (1/2)W_{n+1}+N \kappa^{-2}  \bar u
+N\hat u+N \bar K .
$$
This allows us to finish the proof as that
of Lemma \ref{lemma 8.1.1}.
The lemma is proved.

{\bf Proof of Theorem \ref{theorem 8.4.1}}.    Take   $x\in  B_{\rho_{0}} $
and $r\leq \rho_{0}$. Then
$B_{r}(x)\subset B_{2\rho_{0}}(x) \subset B_{\rho} $ and $2\rho_{0}\leq 1/\nu$. Hence, by taking $x$ as a new origin
and setting   $\rho_{1}=2\rho_{0}$,
from Lemma \ref{lemma 8.2.1} we infer that
$$
r^{\mu}\dashnorm D^{2}u\|_{L_{p}(B_{r}(x))}
\leq  I.
$$
Since $x\in B_{\rho_{0}}$ and $r\leq \rho_{0}$, we also have
$$
r^{\mu}\dashnorm I_{B_{\rho_{0}}}D^{2}u\|_{L_{p}(B_{r}(x))}
\leq  NI.
$$
This inequality is trivially extended for $r\in [\rho_{0},2\rho_{0}]$ and this proves the theorem.

This theorem takes care of interior estimates
in the case of bounded $\gb$ when $p$ can be any number in $(d_{0},\infty)$.
To treat the  case of $\gb$ with rather  poor 
summability properties
we need some preparation in which we restrict
the range of $p$ from above.
The following fact is crucial.

\begin{lemma}
                                \label{lemma 5.22.3}
Let $1<p<q< d$  and $R\in(0,\infty)$.
Set $\mu=q/p$ and $q'=pq/(q-p)$.
Then for any $p'\in[1,q']$ and $u\in E^{1}_{p,\mu}(B_{R})$, we have
\begin{equation}
                                \label{5.22.5}
\| u \|_{E_{p',\mu-1}(B_{R})}
\leq N 
 \|Du\|_{E _{p,\mu}(B_{R})}+
N R^{-1}
 \|u\|_{E _{p,\mu}(B_{R})},
\end{equation}
where the constants $N$ depend only on $d,p,q$.
\end{lemma}

Proof. In light of H\"older's inequality
we may assume that $p'=q'$. Scalings show that we may assume that $R=1$. In that case consider the mapping $\Phi:\bar B_{3/2}\to\bar B_{1}$, $\Phi(x)=x(2/(|x|\vee 1)-1)$ that   preserves $B_{1}$,  
is Lipschitz continuous
and has Lipschitz continuous inverse if restricted to $\bar B_{3/2}\setminus B_{1}$.
  Then, obviously, for any 
$v\in E _{p,\mu}(B_{1})$
\begin{equation}
                                \label{5.22.4}
\|w 
\|_{E_{p,\mu}(B_{3/2})}
\leq N \|v 
\|_{E _{p,\mu}(B_{1})},
\end{equation}
where $N=N(d,p,q)$ and $w(x)=v(\Phi (x))$.

Now take $x\in B_{1}$, $\rho\leq 2$,
and take $\zeta\in C^{\infty}_{0}(\bR^{d})$
such that $\zeta=1$ on $B_{1}$, $\zeta=0$
outside $B_{3/2}$, and $|\zeta|+|D\zeta| 
\leq N=N(d)$.

Since $\mu-1=\mu p/q'$ and $1/q'=1/p-1/(\mu p)$, 
by the Adams theorem 3.1 of \cite{Ad_75} or
Theorem 2.1 of \cite{DK_21} for $w= u (\Phi )$ we have  
$$
\rho^{\mu-1}\dashnorm uI_{B_{1}}\|_{L_{q'}( B_{  \rho}(x))}
\leq N\rho^{\mu-1}\dashnorm \zeta w\|_{L_{q'}( B_{ \rho}(x))}
$$
$$
\leq N\|D (\zeta w)\|_{E_{p,\mu}}
\leq N\|w\|_{E^{1}_{p,\mu}(B_{3/2})}.
$$
It only remains to note that the last expression
is less than the right-hand side of 
\eqref{5.22.5}
in light of \eqref{5.22.4}. The lemma is proved.

\begin{corollary}
                            \label{corollary 5.22.2}
Let $1<p<q< d$, $R\in(0,\infty)$,  and $b\in E_{q,1}(B_{ R})$.
Set $\mu=q/p$.
Then for any $u\in E^{2}_{p,\mu}(B_{R})$, we have
$$
\|b|D u|\|_{E_{p,\mu}(B_{R})}
\leq N\|b\|_{E_{q,1}(B_{R})}\big(
 \|D^{2}u\|_{E _{p,\mu}(B_{R})}
+R^{-1}\|D u\|_{E _{p,\mu}(B_{R})}\big)
$$
\begin{equation}
                              \label{8.5.3}
\leq N\|b\|_{E_{q,1}(B_{R})}\big(
 \|D^{2}u\|_{E _{p,\mu}(B_{R})}
+R^{-2}\| u\|_{E _{p,\mu}(B_{R})}\big),
\end{equation}
where the constants $N$ depend only on $d,p,\mu$.
\end{corollary}

Indeed, by H\"older's inequality 
$$
\rho^{\mu}
\dashnorm I_{B_{R}}b|Du| \|_{L_{p}(B_{ \rho}(x))}
\leq  \rho\dashnorm I_{B_{R}}b \|_{L_{q}(B_{ \rho}(x))}
\rho^{\mu-1}\dashnorm I_{B_{R}}Du\|_{L_{q'}(B_{ \rho}(x))},
$$
where $q'=pq/(q-p)$.
This and \eqref{5.22.5}
obviously lead to the first inequality.
The second one follows from
\eqref{8.4.3}.

{\bf Proof of Theorem \ref{theorem 8.5.1}}. Set   $\rho_{1}=2
\rho_{0}$ and use $\kappa$, $r_{n}$ from the proof
of Lemma \ref{lemma 8.1.1}.
By Theorem \ref{theorem 8.4.1}
with $F[u]+K$ in place of $H[u]$,
where $K=-F[\zeta_{n}u]$, we have 
\begin{equation}
                              \label{8.5.1}
I_{n}:=\|D^{2}u\|_{E _{p,\mu }(B_{r_{n}})}
\leq \|D^{2}(\zeta_{n} u)\|_{E _{p,\mu }(B_{r_{n+1}})}
\leq
J+N_{1}
\| F[\zeta_{n} u]\|_{E_{p,\mu}(B_{\rho_{n+1}})} ,
\end{equation}
where $J=N\rho^{-2} 
\| u\|_{E_{p,\mu}(B_{\rho})}+N  
\rho^{\mu-2}\sup_{B_{\rho}}|u|$.
In light of the arguments in the proof of Lemma \ref{lemma 8.1.1} and \eqref{8.4.3}, \eqref{8.5.2},  
and \eqref{8.5.3},
  the
  last term  in \eqref{8.5.1} is dominated by
$$
\| F[ u]\|_{E_{p,\mu}(B_{r_{n+1}})}
+N 2^{n}\kappa^{-1}
\| Du\|_{E_{p,\mu}(B_{r_{n+1}})}
$$
$$
+N4^{n}\kappa^{-2}\|  u\|_{E_{p,\mu}(B_{\rho_{1}})}
+N\| K\|_{E_{p,\mu}(B_{\rho_{1}})}
$$
$$
\leq N_{1}\big(\hat\theta +
 \hat b+2^{n}\kappa^{-1}\varepsilon\big)
I_{n+1}
$$
$$
+N\big(\hat b r_{n+1}^{-2}+ 4^{n}\kappa^{-2}+2^{n}\kappa^{-1}\varepsilon^{-1}\big)\|  u\|_{E_{p,\mu}(B_{\rho_{1}})}+N\| K\|_{E_{p,\mu}(B_{\rho_{1}})}
$$
for any $\varepsilon\in(0,1]$. We now
choose and fix $\hat\theta$ and $\hat b$
so that $N_{1}(\hat\theta+ \hat b)\leq 1/16$
and take the largest $\varepsilon\in(0,1]$
for which $N_{1}2^{n}\kappa^{-1}\varepsilon\leq 1/16$. Then
$$
\hat b r_{n+1}^{-2}+ 4^{n}\kappa^{-2}+2^{n}\kappa^{-1}\varepsilon^{-1}\leq N4^{n}\kappa^{-2},
$$
so that coming back to \eqref{8.5.1} we get
$$
I_{n}\leq (1/8)I_{n+1}+N4^{n}\kappa^{-2}
\|  u\|_{E_{p,\mu}(B_{\rho_{1}})}+N\| K\|_{E_{p,\mu}(B_{\rho_{1}})}.
$$
This allows us to finish the proof 
as that of Lemma \ref{lemma 8.1.1}.
The theorem is proved.

\mysection{Existence theorems}
                         \label{section 8.11.4}

Here is a general result in which
only a few of our assumptions are supposed to hold.

\begin{lemma}
                         \label{lemma 8.7.2}
Let   $u\in W^{2}_{d_{0},\loc}(\Omega)\cap C(\bar \Omega)$ satisfy \eqref{7.29.1} in $\Omega$     and be such  
that $u=g $ on $\partial\Omega$.
Let estimate \eqref{3.26.20} hold  with $q=d_{0}$ and
$\hat b=\hat b(d,\delta)$ from Theorem 1.1
of \cite{Kr_21a}
for any $r\in(0,R_{0}]$ and ball
$B $ of radius $r$ and let Assumptions 
\ref{assumption 3.11.2} and \ref{assumption 3.11.3} be satisfied.  
Suppose that Assumption  \ref{assumption 10.5.1}    
is satisfied with some $\hat\theta$.
For $\rho>0$ introduce constants $M_{\rho}$
such that
$$
M_{\rho}\geq \|u\|_{W^{2}_{d_{0}}(\Omega^{\rho})}.
$$
Then the modulus of continuity of $u$
in $\bar \Omega$ is dominated by a 
continuous function $\omega(r)$, $r\geq0$, such that $\omega(0)=0$,
depending only on $r$, $R_{0}$, $p$, $M_{\cdot}$, the diameter
of $\Omega$, $\rho,\gamma$ from Assumption
\ref{assumption 3.11.2}, $L_{d_{0}}(\Omega)$-norms of $K$ and $u$, and the modulus
of continuity of $g$.
\end{lemma}  

Proof. Looking at 
$0= \big[H(u,Du,D^{2}u)-H(u,Du,0)\big] 
+H(u,Du,0)$ and using
Assumptions \ref{assumption 3.11.3} and \ref{assumption 10.5.1} one sees
that $a^{ij}D_{ij}u+b^{i}D_{i}u+f=0$,
where $(a^{ij})$ is $\bS_{\delta}$-valued, $|b|\leq \gb $, and $|f|\leq K_{0}|u|+K$.
Then $|u(x_{1})-u(x_{2})|$
for $x_{1},x_{2}$ that are close to
$\partial \Omega$ is estimated by using Lemma
\ref{lemma 7.19.1}. If they are far,
we use embedding theorems ($d_{0}>d/2$) to estimate the
difference in terms of $M_{\cdot}$.
The combination of these estimates
leading to the desired result is a simple
exercise. The lemma is proved.

Coming closer to the proof of Theorem
\ref{theorem 10.5.1},
observe that estimate \eqref{1.8.2} follows
from Theorem \ref{theorem 4.1.1}.
Indeed, the assumption of this theorem
concerning \eqref{3.26.20} is obviously satisfied since $\gb\leq K_{0}$ and on
 the set $\Omega\cap\{u>0\}$ we have
$$
0=H[u]=\big[H(u,Du,D^{2}u)-H(u,Du,0)\big]
+H(u,Du,0)
$$
$$
\leq a^{ij}D_{ij}u+\gb|Du|+K=
a^{ij}D_{ij}u+b^{i} D_{i}u +K,
$$
where $(a^{ij})$ is $\bS_{\delta}$-valued and $|b|\leq \gb\leq K_{0}$. Hence,
$$
u\leq N\|I_{\Omega,u>0}K\|_{L_{d_{0}} }
+\sup_{\partial\Omega}u_{+}.
$$
Similarly, the estimate of $-u$ is obtained.

In the following lemma $\hat H(\sfu'',x)$
is a measurable function such that it is Lipschitz
continuous with respect to $\sfu''$,
$(D_{\sfu''_{ij}}\hat H)\in\bS_{\delta}$ at all
points of differentiability of $\hat H$
and $\hat H[0]\in L_{p,\loc}(\bR^{d})$.
Observe that $\hat H(\sfu'',\cdot)\in L_{p,\loc}(\bR^{d})$
for any $\sfu''$ and by Lebesgue's
theorem
\begin{equation}
                        \label{8.6.5}
\lim_{r\downarrow 0}\dashnorm
\hat H(\sfu'',\cdot)-
\hat H(\sfu'',x_{0})\|_{L_{p}(B_{r}(x_{0}))}=0
\end{equation}
for almost any $x_{0}$. Since $\hat H(\sfu'',x)$
is Lipschitz continuous in $\sfu''$ one can choose a set of $x_{0}$ of full measure
such that \eqref{8.6.5} holds for any $\sfu''$.
\begin{lemma}
                          \label{lemma 8.6.1}
 
Let $p\in[d_{0},\infty)$, $R\in(0,\infty)$,
$u\in W^{2}_{p}(B_{R})$, $f\in L_{p}(B_{R})$ . Then $\hat H[u]\geq f$  in $B_{R}$ if and only if
for any $B_{r}(x_{0})\subset B_{R}$ and any
$\phi\in C^{2}(B_{r}(x_{0}))$ we have in $B_{r}(x_{0})$ that
\begin{equation}
                              \label{8.6.3}
u\leq \phi+
\sup_{\partial B_{r}(x_{0})}(u-\phi)_{+}
+Nr^{2}\dashnorm (f-\hat H[\phi])_{-}\|_{L_{p}(B_{r}(x_{0}))}
\end{equation}
with $N$ independent of $u,\phi,x_{0},r$.  

\end{lemma}

Proof. ``only if''. We have $f-\hat H[\phi]
\leq H[u]-\hat H[\phi]=a^{ij}D_{ij}(u-\phi)$ with a $\bS_{\delta}$-valued $(a^{ij})$ and \eqref{8.6.3}
follows from Theorem \ref{theorem 4.1.1}
and a scaling argument.

``if''.  Take   $x_{0}\in B_{R}$ such that
$$
u(x)=u(x_{0})+(x-x_{0})^{i}D_{i}u(x_{0})
+(1/2)(x-x_{0})^{i}(x-x_{0})^{j}D_{ij}u(x_{0})
+o(|x-x_{0}|^{2}).
$$
By the  Zygmund-Calder\'on theorem one
can take almost any $x_{0}\in B_{R}$ since $p>d/2$. We may even restrict this $x_{0}$
to satisfy \eqref{8.6.5} for any $\sfu''$
and satisfy
$$
\lim_{r\downarrow 0}\dashnorm f-f(x_{0})\|
_{L_{p}(B_{r}(x_{0}))}=0.
$$
Then fix an $\varepsilon>0$ and for
$r$, such that $o(r^{2})\leq \varepsilon
r^{2}$, in $B_{r}(x_{0})$ introduce 
$$
\phi(x)=u(x_{0})+(x-x_{0})^{i}D_{i}u(x_{0})
$$
$$
+(1/2)(x-x_{0})^{i}(x-x_{0})^{j}D_{ij}u(x_{0})
+\varepsilon\big(2|x-x_{0}|^{2}-r^{2}\big).
$$
We will send $r\downarrow 0$ and, therefore,
we may concentrate on $r$ such that
$B_{r}(x_{0})\subset B_{R}$ and 
$u\leq \phi$ on $\partial B_{r}(x_{0})$.
Then \eqref{8.6.3} at $x=x_{0}$ yields
$$
\varepsilon r^{2}\leq Nr^{2}\dashnorm\big(
f-\hat H(D_{ij}u(x_{0})+4\varepsilon \delta^{ij} ,\cdot)\big)_{-}\|_{L_{p}(B_{r}(x_{0}))},
$$
which after letting $r\downarrow0$ becomes
$$
\varepsilon\leq N\big(
f(x_{0})-\hat H(D_{ij}u(x_{0})+4\varepsilon \delta^{ij} ,x_{0})\big)_{-}.
$$
We send $\varepsilon\downarrow0$ and get 
$f(x_{0})-\hat H(D_{ij}u(x_{0}) ,x_{0})\leq 0$, thus proving the lemma.

Similarly, or just taking $-\hat H(-\sfu'',x)$
in place of $\hat H(\sfu'',x)$, one proves the following.

\begin{lemma}
                          \label{lemma 8.7.1}
 
Let $p\in[d_{0},\infty)$, $R\in(0,\infty)$,
$u\in W^{2}_{p}(B_{R})$, $f\in L_{p}(B_{R})$. Then $\hat H[u]\leq f$ in $B_{R}$ if and only if  for any $B_{r}(x_{0})\subset B_{R}$ and any
$\phi\in C^{2}(B_{r}(x_{0}))$ we have in $B_{r}(x_{0})$ that
\begin{equation}
                              \label{8.7.1}
u\geq \phi-
\sup_{\partial B_{r}(x_{0})}(u-\phi)_{-}
-Nr^{2}\dashnorm (f-\hat H[\phi])_{+}\|_{L_{p}(B_{r}(x_{0}))}
\end{equation}
with $N$ independent of $u,\phi,x_{0},r$.

\end{lemma}

{\bf Proof of Theorem \ref{theorem 10.5.1}}.
For $n=1,2,...$ introduce
$$
H_{n}(\sfu,x)=I_{K\leq n}H (\sfu  ,x)+
I_{K>n}F(\sfu'',x).
$$
Observe that $H_{n}=F+G_{n}$, where
$
G_{n}(\sfu,x)=I_{K\leq n}G(\sfu ,x),
$
$$
|G_{n}(\sfu,x)|\leq \hat\theta|\sfu''|+K_{0}
|\sfu'  |+ 
I_{K\leq n}K.
$$
Here   the free term belongs to $L_{p,\loc}(\bR^{d})$ for any $p>1$. It follows from  
10.1.14 of \cite{Kr_18} (see also
Remark 10.1.15 there) that for appropriate 
$\hat\theta, \theta$, 
 depending only on
$d,\delta,p$, under Assumptions
\ref{assump1} and \ref{assumption 10.5.1},
for any $n$, there exists $u\in W^{2}_{p+d,\loc}(\Omega)\cap C(\bar\Omega)$ such that
$H_{n}[u_{n}]=0$ in $\Omega$ and $u_{n}=g$
on $\partial\Omega$. 

Estimate \eqref{1.8.2}
shows that $u_{n}$ are uniformly bounded in
$\bar\Omega$. Then Theorem \ref{theorem 8.4.1}
implies that for any $\rho>0$ the 
$E^{2}_{p,\mu}(\Omega^{\rho})$-norms
of $u_{n}$ are bounded, provided that
  Assumptions
\ref{assump1} and \ref{assumption 10.5.1}
are satisfied with appropriate $\hat\theta,\theta$. In particular, by Lemma \ref{lemma 8.7.2}
the family $\{u_{n}\}$ is uniformly bounded
and uniformly continuous in $\bar\Omega$.
Therefore, there exists a subsequence
$u_{n_{k}}$ and $u\in C(\bar\Omega)$ such that
$u_{n_{k}}\to u$ uniformly on $\bar \Omega$.

Next, by the compactness of embeddings,
for each $\rho>0$, the family $Du_{n}$
is precompact in $L_{p}(\Omega^{\rho})$.
Hence, by using the 
Cantor diagonalization method
we may assume that  
\begin{equation}
                             \label{8.7.4}
\tilde K:= |Du_{n_{1}}|+\sum_{k=1}^{\infty} |Du_{n_{k+1}}-Du_{n_{k}}|
\in L_{p}(\Omega^{\rho})
\end{equation}
for any $\rho>0$. Then $Du_{n_{k}}$ converge in
$L_{p}(\Omega^{\rho})$ for any $\rho>0$ and
almost everywhere in $\Omega$ to some functions
which automatically coincide with $Du$.
The weak limit  of $D^{2}u_{n_{k}}$ is, of course, $D^{2}u$, so that $u\in E^{2}_{p,\mu,\loc}(\Omega)$.
 
To prove that $H[u]=0$ in $\Omega$, for
$m=1,2,...$ set
$$
\hat H_{m}(\sfu'',x)=\sup_{k\geq m}H_{n_{k}}
(u_{n_{k}}(x),Du_{n_{k}}(x),\sfu'',x).
$$
Observe that
\begin{equation}
                             \label{8.8.2}
|\hat H_{m}(0,x)|\leq K_{0}\sup_{n=1,2,...}
\sup_{\Omega}|u_{n}|+ \tilde K+K\in L_{p,\loc}(\Omega).
\end{equation}
Also, obviously, $(D_{\sfu''_{ij}}\hat H_{m})
\in \bS_{\delta}$ at all points of differentiability of $\hat H_{m}$.
For $k\geq m$ we have $\hat H_{m}(D^{2}u_{n_{k}})\geq0$ in $\Omega$ implying by Lemma
\ref{lemma 8.6.1} that for any $\bar B_{r}
(x_{0})\subset \Omega$ and any
$\phi\in C^{2}(B_{r}(x_{0}))$ we have in $B_{r}(x_{0})$ that
\begin{equation}
                              \label{8.8.1}
u_{n_{k}}\leq \phi+
\sup_{\partial B_{r}(x_{0})}(u_{n_{k}}-\phi)_{+}
+Nr^{2}\dashnorm ( \hat H_{m}[\phi])_{-}\|_{L_{p}(B_{r}(x_{0}))}
\end{equation}
with $N$ independent of $u,\phi,x_{0},r$,
and $m$. The fact that $N$ is indeed independent of $m$ easily follows from the proof of Lemma \ref{lemma 8.6.1}. We   pass to the limit as $k\to\infty$, which allows us to replace $u_{n_{k}}$ on the left in \eqref{8.8.1} with $u$ and conclude by Lemma \ref{lemma 8.6.1} that $\hat H_{m}(D^{2}u)\geq0$ in
$\Omega$. This inequality on the set
$\{K\leq n_{m}\}$ means that
$$
\sup_{k\geq m}H 
(u_{n_{k}}(x),Du_{n_{k}}(x),D^{2}u(x),x)\geq0.
$$ 
Since $Du_{n_{k}}(x)\to Du$ 
almost everywhere in $\Omega$, $u_{n_{k}}(x)\to u(x)$ in $\Omega$ and $H$ is a continuous function of $\sfu'$, by setting $m\to\infty$, we conclude
that $H[u]\geq0$ in $\Omega$. Similarly,
by using Lemma \ref{lemma 8.7.1} one proves that $H[u]\leq0$ in $\Omega$. The theorem is proved.

{\bf Proof of Theorem \ref{theorem 8.8.1}}. Estimate \eqref{8.8.3} is derived as
(the identical) \eqref{1.8.2} by using
Theorem \ref{theorem 4.1.1}. To prove the
existence,
for $n=1,2,...$ introduce
$$
H_{n}(\sfu,x)= H (\sfu'_{0},
n[\sfu']/(n+\gb),\sfu'',x).
$$
Observe that $H_{n}=F+G_{n}$, where
$$
|G_{n}(\sfu,x)|\leq \hat\theta|\sfu''|+K_{0}
|\sfu'_{0}|+\big(n\gb/(n+\gb)\big)|[\sfu']|+
 K.
$$
We apply Theorem \ref{theorem 10.5.1} upon observing that above  the coefficients of $\sfu'$ are bounded
and the free term belongs to 
$E_{p,\mu,\loc}(\Omega)\cap L_{p}(\Omega)$.
Then we conclude that there exists
$u_{n}\in E^{ 2}_{p,\loc}(\Omega)\cap C(\bar\Omega)$
 satisfying \eqref{7.29.1} in $\Omega$ 
with $H_{n}$ in place of $H$ and such  
that $u_{n}=g $ on $\partial\Omega$.
Theorem \ref{theorem 8.5.1} guarantees
that for any $\rho>0$ the 
$E^{2}_{p,\mu}(\Omega^{\rho})$-norms
of $u_{n}$ are bounded, provided that
  Assumptions \ref{assumption 10.5.1},
\ref{assump1}, and \ref{assumption 4.29.1}  
are satisfied with appropriate $\hat\theta,\theta,\hat b$. After that, as in the proof of
Theorem \ref{theorem 10.5.1}, we find
a subsequence
$u_{n_{k}}$ and $u\in W^{1}_{p,\loc}(\Omega)\cap C(\bar\Omega)$ such that
$u_{n_{k}}\to u$ uniformly on $\bar \Omega$
and $Du_{n_{k}}\to Du$ almost everywhere.
Of course, $u\in E^{2}_{p,\mu,\loc}(\Omega)$.

Furthermore, by Corollary \ref{corollary 5.22.2}, for any $B_{r}(x_{0})\subset\Omega$
with $r\leq R_{0}$ the integrals
$$
\int_{B_{r}(x_{0})}\gb^{p}|Du_{n}|^{p}\,dx
$$
are bounded by a constant independent of $n$.
It follows that for any $p'\in[d_{0},p)$
$$
\lim_{k\to\infty}\int_{B_{r}(x_{0})}\gb^{p'}|Du_{n_{k}}-Du|^{p'}\,dx=0
$$
and there exists a subsequence which we
identify with the above one such that
\begin{equation}
                             \label{8.8.4}
\tilde K:=\gb|Du_{n_{1}}|+
\sum_{k=1}^{\infty}\gb|Du_{n_{k+1}}-Du_{n_{k}}|
\in L_{d_{0},\loc}(\Omega ).
\end{equation}
Then we introduce $\hat H_{m}$ as in the proof
of Theorem \ref{theorem 10.5.1}, observe
that \eqref{8.8.2} holds with $d_{0}$
in place of $p$   and with the help
of Lemma \ref{lemma 8.6.1} conclude
that $\hat H_{m}(D^{2}u)\geq0$ in
$\Omega$. Since $u_{n_{k}},n_{k} Du_{n_{k}}/(n_{k}+\gb)\to u, Du$ almost everywhere as $k\to
\infty$ and $H$ is continuous in $\sfu'$,
$$
0\leq \lim_{m\to\infty}\hat H_{m}(D^{2}u)
=\nlimsup_{k\to\infty}
H(u_{n_{k}},n_{k} Du_{n_{k}}/(n_{k}+\gb),D^{2}u)
=H[u].
$$
Similarly,
by using Lemma \ref{lemma 8.7.1} one proves that $H[u]\leq0$ in $\Omega$. The theorem is proved.

\end{document}